\documentclass{amsart}
\usepackage[latin1]{inputenc}
\usepackage[T1]{fontenc}
\usepackage{lmodern}
\usepackage[english]{babel}
\usepackage{microtype}
\usepackage[fontsize=11pt]{scrextend}
\usepackage{amsmath,amssymb,amsfonts,amsthm}
\usepackage{mathtools,accents}
\usepackage{mathrsfs}
\usepackage{aliascnt}
\usepackage{braket}
\usepackage{bm}
\usepackage[initials]{amsrefs}
\usepackage[citecolor=blue,colorlinks]{hyperref}
\addto\extrasenglish{}

\usepackage{enumerate}
\usepackage{xcolor}

\usepackage{aliascnt}

\makeatletter
\def\newaliasedtheorem#1[#2]#3{
  \newaliascnt{#1@alt}{#2}
  \newtheorem{#1}[#1@alt]{#3}
  \expandafter\newcommand\csname #1@altname\endcsname{#3}
}
\makeatother

\theoremstyle{plain}
\newtheorem{theorem}{Theorem}[section]
\newaliasedtheorem{lemma}[theorem]{Lemma}
\newtheorem{prop}[theorem]{Proposition}
\newtheorem{corollary}[theorem]{Corollary}
\newtheorem{remark}[theorem]{Remark}

\numberwithin{equation}{section}

\def\R{\mathbb R}
\def\C{{\mathbb C}}

\DeclareMathOperator{\RE}{Re}
\DeclareMathOperator{\IM}{Im}

\title[Blow-up 3D Davey-Stewartson]{On finite time blow-up for a 3D Davey-Stewartson system} 

\author[L. Forcella]{Luigi Forcella}

\address{Luigi Forcella\medskip\hfill\break  Heriot-Watt University, Department of Mathematics,
Edinburgh, EH14 4AS, United Kingdom,
and The Maxwell Institute for the Mathematical Sciences\medskip}
\email{l.forcella@hw.ac.uk}

\subjclass[2000]{35Q55, 35B40, 35B44}

\keywords{Davey-Stewartson system, NLS-type equation, blow-up}

\begin{document}

\maketitle  


\begin{abstract}
We consider the elliptic-elliptic Davey-Stewartson system in the three-dimensional Euclidean space, and we give sufficient conditions for the existence of finite time blow-up solutions in non-isotropic  spaces. The proof is based on some general results on distributions defined via homogeneous symbols, in conjunction with a convexity argument. 
\end{abstract}

\section{Introduction}
In this short note, we consider the following initial value problem: 
\begin{equation}\label{DS}
\left\{ 
\begin{aligned}
i\partial_{t}u+\Delta u+c_1|u|^{\alpha}u+c_2  E_1(|u|^2)u=0\\
u(0,x)=u_0(x) \in H^1(\R^3)
\end{aligned}\right.,
\end{equation}
where $(t,x)\in [0,T)\times \R^3$, $u:[0,T)\times \R^3\mapsto \C$, $c_1$ and $c_2$ are two positive parameters, $\alpha\in(0,4)$, and the operator $  E_1$ is given in term of the Fourier symbol $\sigma_1(\xi)=\frac{\xi^2_1}{|\xi|^2}$, with $\xi=(\xi_1,\xi_2,\xi_3)\in\R^3$:
\begin{equation}\label{eq:sigma1}
 E_1 f(x)=\mathcal F^{-1}\left(\sigma_1(\cdot)\hat f(\cdot)\right)(x), 
\end{equation}
where $\hat f=\mathcal F f$ denotes the Fourier transform of $f$, and $\mathcal F^{-1}$ stands for the inverse Fourier transform. The equation \eqref{DS} is called  3D Davey-Stewartson system. Though it is a single equation, one refers to it as a system for it can be viewed as a three dimensional extension of the following Davey-Stewartson system, see \cite{DaSt,CSS}: 
\begin{equation}\label{DS-sy}
\left\{ 
\begin{aligned}
i\partial_{t}v+a_1 \partial_{xx}^2v +\partial_{yy}^2v&=a_2 |v|^\alpha v+a_3 v\partial_x w\\
\partial_{xx}^2w+a_4\partial_{yy}^2w&=\partial_{x}(|v|^2)
\end{aligned}\right.,
\end{equation}
where $v=v(t,x,y)$ and $w=w(t,x,y)$, with $(t,x,y)\in\R\times \R\times \R$,  $a_i$, with $i\in\{1,2,3,4\}$, are real parameters. According to the signs of the coefficients $a_1$ and $a_4$, the system \eqref{DS-sy} is classified as: $(+,+)$ elliptic-elliptic, $(+,-)$ elliptic-hyperbolic,  $(-,+)$ hyperbolic-elliptic,  $(-,-)$ hyperbolic-hyperbolic, respectively.  This paper concerns  the elliptic-elliptic case, namely $a_1>0$ and $a_4>0$. See \cite{CSS, DaSt, GhiSa, SS, ZS, NiSa, PSSW}, and references therein for physical insights on the model. From now on, we omit the space $\R^3$, as we work in the 3D case. \medskip

Existence of solutions to the Cauchy problem  \eqref{DS} in the energy space $H^1$ was established in \cite{Guo-Wang}. A solution $u\in C((-T_{-},T_{+});H^1)$ conserves the mass and the energy; specifically, with $u(t)=u(t,x)$, the quantities 
\begin{equation*}
M(u(t)):=\int|u(t)|^2 dx
\end{equation*}
and
\begin{equation*}\label{energy}
E(u(t)):=\frac{1}{2}\int |\nabla u(t)|^2-\frac{c_1}{\alpha+2}\int |u(t)|^{\alpha+2} dx-\frac{c_2}{4}\int   E_1(|u(t)|^2)|u(t)|^2 dx
\end{equation*} 
do not depend on time, for any $t\in(-T_{-},T_{+}),$ where $T_{-},T_{+}\in(0,\infty]$ denote the minimal and maximal time of existence, respectively. In this paper, we are interested in establishing sufficient conditions leading to formation of singularities in finite time for solutions to \eqref{DS}. In particular we prove finite time blow-up for non-isotropic  initial data.   Let us observe that the equation \eqref{DS} does not enjoy a radial invariance; nonetheless, it satisfies a cylindrical invariance -- in the $x_1$-direction -- as the symbol defining $  E_1$ suggests. Hence, we consider the non-isotropic space of cylindrical functions in $H^1$ with a weight in $L^2$. Specifically, by denoting a vector $x\in\R^3$  as $x=(x_1,x_2,x_3)=(x_1,\bar x)$, with $\bar x=(x_2,x_3)$, we introduce the space
\begin{equation*}
\Sigma_1 =\left\{ f \in H^1 : f(x)=f(x_1, |\bar x|) \ \hbox{ and } \  f\in L^2(x_1^2 dx) \right\};
\end{equation*}
$\Sigma_1$ is therefore the sub-space of $H^1$ consisting of functions with radial invariance with respect to the $\bar x$ coordinate, and with finite variance in the $x_1$ direction. As the equation \eqref{DS} is not radial-invariant,  we cannot rely on a radiality assumption to prove finite time blow-up. Indeed, to the best of our knowledge, all the existing literature treating the problem of finite time blow-up deals with solutions with finite variance, i.e. for initial data $u_0$ in $H^1\cap L^2(|x|^2;dx)$. Therefore, our main result is somehow ``minimal'' with respect to the symmetry assumptions on the solutions. Besides the symmetry hypothesis and the finiteness of the variance, sufficient condition for blow-up are given by also imposing some bounds on the initial data. In particular, such conditions are defined in terms of solutions of the elliptic equation associated to \eqref{DS}: 
\begin{equation}\label{eq:ell}
-\Delta Q+Q-c_1|Q|^\alpha Q-c_2  E_1(|Q|^2)Q=0.
\end{equation}
Note that $u(t,x)=e^{it}Q(x)$, where $Q$ solves \eqref{eq:ell}, is a solution to \eqref{DS}, and it is called standing wave solution. Solutions $Q$ to \eqref{eq:ell} allow us to introduce the concept of Ground State. To this aim, we denote the (conserved along the flow) Lagrangian $S(u)=E(u)+\frac12M(u)$, and we denote the set $\mathbb G$ of Ground  States as the set of non-trivial solutions  to \eqref{eq:ell} minimizing the Lagrangian functional:
\[
\mathbb G=\{G\neq0 \hbox{ solving }  \eqref{eq:ell}: S(G)\leq S(Q) \hbox{ for any } Q\neq0 \hbox{ solving }  \eqref{eq:ell}\}.
\]
We refer to \cite{Cipo} for the existence theory of Ground States for \eqref{DS}. Let us observe that a solution $Q$ to \eqref{eq:ell} satisfies $P(Q)=0$, where $P$ is the Pohozaev functional defined by 
\begin{equation}\label{eq:G}
P(f)=\int|\nabla f|^2-\frac{3c_1\alpha}{2(\alpha+2)}\int|f|^{\alpha+2} dx-\frac{3c_2}{4}\int   E_1(|f|^2)|f|^2 dx.
\end{equation} 

With all the notions above, we are able to give our main result. 
\begin{theorem}\label{thm:main} Let $\alpha\in[4/3,2]$ and let $u_0\in\Sigma_1$. Assume that $S(u_0)< S(G)$ where $G$ is a Ground State for \eqref{eq:ell}, and that $P(u_0)<0$. Then the solution to \eqref{DS} blows-up in finite time, i.e. $T_{-}$ and $T_{+}$ are both finite. 
\end{theorem}

In order to prove our main result, we rely on some decay properties for an operator $T$ defined by means of homogeneous symbols of order zero. Specifically, when we pair, in $L^2$, a function $g$ with $Tf$, where the supports of $f$ and $g$ are disjoint with a positive distance $\delta$, then we get a decay of order $\delta^{-3}$. This fact, jointly with a careful localization of the solutions, allows us to employ a convexity argument. \medskip

For the localization argument, we are inspired by the work of Lu and Wu, see \cite{LW}. In their paper, the authors prove scattering results for the system \eqref{DS} (via a concentration/compactness and rigidity scheme), and they state a grow-up result: in particular, they show that provided the initial datum $u_0$ satisfies only $S(u_0)< S(G)$ and $P(u_0)<0$, if the solution to \eqref{DS} is global, then there exists  a diverging sequence of times $\{t_n\}_n$ such that $\limsup_{n\to\infty}\|u(t_n)\|_{\dot H^1}=\infty$.  But they only prove finite time blow-up for solution in $H^1\cap L^2(|x|^2dx)$. 
Moreover, their result is based on point-wise type decay estimates for $E_1$. We prove our result by employing  much simpler estimates coming from general properties of operators with homogeneous symbols, see Proposition \ref{thm:zero-degree}.\medskip

We give some remarks. 
\begin{remark}
\textup{As already mentioned, the fact that $E_1$  does not leave invariant the set of radial functions prevents us to give a result for radial initial data. Due to the structure of the symbol $\sigma_1$ defining $E_1$, see \eqref{eq:sigma1}, we can instead show the finite time blow-up result for cylindrical solutions. Therefore, our hypothesis is not purely artificial, and as pointed-out in \cite[Remark 3.8]{Cipo}, it is also linked to the possible existence of Ground States with such a symmetry. }
\end{remark}
\begin{remark}
\textup{We cover the range of non-linearities $\alpha\in[4/3,2]$. The lower bound corresponds to the mass-critical case. The upper bound instead plays the same role of  the limitation $\alpha\leq4$ for radial solutions for the 2D NLS equation, see Ogawa and Tsutsumi \cite{OgTs}. }
\end{remark}
\begin{remark}
\textup{To the best of our knowledge, all the results on formation of singularities in finite time concern  solutions with finite variance, see \cite{GanZha, GhiSa, LW, LLY, Lu2016, ZZ2019, Zhu, ZZ2011}. }
\end{remark}
\begin{remark}
\textup{A large amount of works have been devoted to the existence and  dynamics of solutions for Davey-Stewartson systems, both in 2D and 3D: we refer the readers to  \cite{GhiSa, Oza, Cipo, Cipo93, Ohta1, Ohta2} and references therein. }
\end{remark}

\section{Preliminary tools}

As mentioned in the Introduction, we will employ a convexity  argument to prove our main result. To this aim, we strongly rely on the following general result, which will enable us treat to the non-local terms -- coming from the presence of non-local non-linearity in the equation -- in the virial estimates. \medskip

We consider a pseudo-differential operator $T$ defined by means of a symbol $\sigma(\xi)$, i.e. $Tf=\mathcal F^{-1}(\sigma \hat f)$, where $\sigma$ is homogenous of order zero, namely $\sigma(\lambda \xi)=\sigma(\xi)$ for any $\lambda>0$, and it is smooth in $\R^3_\xi\setminus\{0\}$. Hereafter, $\langle\cdot,\cdot\rangle$ will denote the $L^2$ pairing. We have the following. 
\begin{prop}\label{thm:zero-degree}
Let $T$ defined as above. Let $f,g\in L^1$ have disjoint supports, and suppose that $\gamma:=\hbox{distance}(\hbox{supp }(f), \hbox{supp }(g))>0$. Then 
\begin{equation}\label{eq:int-decay-general}
|\langle Tf,g\rangle|\lesssim \gamma^{-3}\|g\|_{L^1}\|f\|_{L^1}.
\end{equation}
\end{prop}

\noindent Observe that the symbol $\sigma_1$ defining $E_1$  fulfils the hypothesis of Proposition \ref{thm:zero-degree}. Moreover, the operator $E_1^2$ defined by means of the symbol $\sigma_1^2(\xi)=\frac{\xi_1^4}{|\xi|^4}$ fall down into the same scenario. Therefore we have the following corollary for functions supported on disjoint cylinders.  
\begin{corollary}\label{cor:decay}
Let $f,g\in L^1$ such that   $\hbox{supp }(f)\subset \{|\bar x|\geq \gamma_2R\}$ and $\hbox{supp }(g)\subset\{|\bar x|\leq \gamma_1R\},$ where $\gamma_1$ and $\gamma_2$ are positive parameters satisfying $\gamma_2-\gamma_1>0.$ Then, for $k=1,2$,
\begin{equation}\label{eq:est-r4}
|\langle E_1^{k}f,g\rangle|\lesssim R^{-3}\|g\|_{L^1}\|f\|_{L^1}.
\end{equation}
\end{corollary}
The proof of the above result was given by Bellazzini and the author in \cite{BF-CV}, where we studied another NLS-type equation with non-local nonlinearity. As we would like to keep this note self-contained, we report the proof for sake of completeness. 
\begin{proof} Under the structural hypothesis for $T$, we have by \cite[Proposition 2.4.7]{Grafakos} that there exist a smooth function $\Omega$ on the two-dimensional sphere $\{z\in\R^3 : |z|=1\}$, and a complex number  $c$ such that 
\[
(\mathcal F^{-1}\sigma)(x) =\frac{1}{|x|^3}\Omega\left(\frac{x}{|x|}\right)+c\delta(x),
\]  
where $\delta$ is the Dirac delta. 
Hence, 
\[
\begin{aligned}
\langle Tf,g\rangle&= \int \left(\frac{1}{|\cdot |^3}\Omega\left(\frac{\cdot}{|\cdot|}\right)\ast f\right)(x)\bar g(x) dx+c\int \left(\delta \ast f\right)(x)\bar g(x) dx\\
&=\iint \frac{1}{|x-y|^3}\Omega\left(\frac{x-y}{|x-y|}\right)f(y)\bar g(x)  dy dx,
\end{aligned}
\]
where the term with the Dirac delta disappears due to the disjointness of the supports. Therefore, as $|x-y|\geq|\bar x|-|\bar y|\geq (\gamma_2-\gamma_1)R$, we have 
\[
|\langle Tf,g\rangle|\lesssim R^{-3}\|\Omega\|_{L^\infty}\|f\|_{L^1}\|g\|_{L^1}\lesssim  R^{-3}\|f\|_{L^1}\|g\|_{L^1}.
\]
\end{proof} 

\begin{remark}\textup{The  general result above allows to avoid point-wise estimates as in \cite{BF19, BF-CV, DFH, LW}, hence simplifying the proofs in the latter papers. 
}
\end{remark}

The next two Propositions are contained in \cite{LW}, in particular see \cite[Corollary 2.7]{LW} and \cite[Corollary 2.9]{LW}, respectively. They are consequences of the variational characterization of the Ground States. 

\begin{prop}\label{prop:P-nega} Let $u_0$ be an initial datum satisfying $S(u_0)<S(G)$ and $P(u_0)<0$. Then  the corresponding solution $u(t)$ to \eqref{DS} satisfies the same bounds, namely $S(u(t))<S(G)$ and $P(u(t))<0$ for any $t\in(-T_{-}, T_{+})$.
\end{prop}
\begin{prop}\label{prop:P-grad} Let $u_0$ be an initial datum satisfying $S(u_0)<S(G)$ and $P(u_0)<0$. Then  there exist $\varepsilon>0$ and $\bar\varepsilon>0$ such that the corresponding solution $u(t)$ to \eqref{DS} satisfies $S(u(t))<(1-\varepsilon)S(G)$ and $P(u(t))<-\bar\varepsilon\|u(t)\|_{\dot H^1}^2$ for any $t\in(-T_{-}, T_{+})$.
\end{prop}
\noindent In particular, the latter Proposition will be crucial when employing a convexity argument to show the blow-up. \medskip

We conclude this section by reporting the following embedding. For any cylindrical function $f\in H^1$ we have 
\begin{equation}\label{eq:strauss} 
\|f\|_{L^4(|\bar x|\gtrsim R)}^4\lesssim R^{-1}\|f\|_{\dot H^1}^2.
\end{equation}
A proof can be found in \cite{BF-CV}, and it is based on the Strauss embedding for radial functions. 
\begin{remark}\label{rem:bound}\textup{ It is worth mentioning that $E_1^k$, $k=1,2$, are $L^2\mapsto L^2$ continuous operators. The latter property will be often used during the rest of the paper, and it easily follows by the boundedness of their symbols. More in general, they are $L^p\mapsto L^p$ continuous for any $p\in(1,\infty)$, see \cite{Cipo}. 
}
\end{remark}

\section{Proof of main result}
In this Section, we give a proof of Theorem \ref{thm:main}. It will be done by using the technical tools introduced in the previous Section, and it relies on a virial argument, together with appropriate localizations of the solution, which enable us to use  Proposition \ref{thm:zero-degree} to control various non-local terms. \medskip

Given a smooth, non-negative, real function $\rho=\rho(x)$ defined on $\R^3$, we define, for a solution $u=u(t,x)$ to \eqref{DS} (we omit the space-time dependence), the time depending function
\begin{equation*}
V_{\rho}(t)=\int \rho|u|^2 dx.
\end{equation*} 
Usual computations, which may be justified by a regularization argument, yield 
\begin{equation}\label{virial1}
\begin{aligned}
\frac{d}{dt}V_{\rho}(t)&=2\IM\left\{\int\nabla \rho\cdot\nabla u\bar u dx\right\},
\end{aligned}
\end{equation}
where we used the equation satisfied by $u$.
By using \eqref{virial1} and again the equation solved by $u$, we have 
\begin{equation*}
\begin{aligned}
\frac{d^2}{dt^2}V_{\rho}(t)&=4 \RE \int\left(\nabla^2\rho\cdot\nabla u\right)\cdot\nabla\bar u dx-\int\Delta^2\rho|u|^2 dx\\
&-\frac{2c_1\alpha}{\alpha +2}\int\Delta \rho|u|^{\alpha+2} dx+2c_2\int\nabla \rho\cdot\nabla\left(  E_1(|u|^2)\right)|u|^2 dx.
\end{aligned}
\end{equation*}

\noindent We precisely chose a function $\rho$ to fit with our symmetry assumptions on the solution.  We consider $\psi: \R^2 \to \R$  a smooth radial function, and by setting $\rho(x) = x_1^2+\psi_R(\bar x)$, with the resclaling $\psi_{R}=R^2\psi(|\bar x|^2/R^2)$, since $u(t) \in \Sigma_1 $ for all $t\in (-T_-,T_+)$,  we have
\begin{align*}
\frac{d^2}{dt^2}V_{x_1^2+\psi_R(\bar x)}(t)&= -\int \Delta^2_{\bar x} \psi_R(\bar x) |u|^2 dx + 4\int \psi_R''(r) |\nabla_{\bar x} u|^2  dx\\
& + 8 \|\partial_{x_1} u\|^2_{L^2} -\frac{2c_1\alpha}{\alpha +2}\int (2+\Delta_{\bar x}\psi_R)|u|^{\alpha+2}dx\\
&+2c_2 \int\nabla_{\bar x} \psi_R\cdot\nabla_{\bar x}\left(  E_1(|u|^2)\right)|u|^2 dx\\
&+4c_2\int x_1\partial_{x_1}\left(  E_1(|u|^2)\right)|u|^2 dx.
\end{align*}
The subscript $\bar x$ above and in what follows means that the differential operator is taken only with respect the $\bar x$ variables. By straightforward  computations, we get 
\begin{align}\label{vir0}
\frac{d^2}{dt^2} V_{x_1^2+\psi_R(\bar x)}(t)&= 8\left(\int |\nabla u|^2-\frac{3c_1\alpha}{2(\alpha+2)}\int |u|^{\alpha+2} dx\right) \\\label{vir-1}
& -\int \Delta^2_{\bar x} \psi_R |u|^2 dx- 4\int(2- \psi_R''(r)) |\nabla_{\bar x} u|^2  dx\\\label{vir-2}
&+\frac{2c_1\alpha}{\alpha +2}\int (4-\Delta_{\bar x}\psi_R)|u|^{\alpha+2}dx\\\label{vir-3}
&+ 2c_2\int\nabla_{\bar x} \psi_R\cdot\nabla_{\bar x}\left(  E_1(|u|^2)\right)|u|^2 dx\\\label{vir-4}
&+4c_2\int x_1\partial_{x_1}\left(  E_1(|u|^2)\right)|u|^2 dx.
\end{align}

\noindent 	By following  Martel \cite{Mar}, we define 
\[
\psi(r)=r-\int_0^r(r-s)\eta(s)  ds,
\]
where the real regular function $\eta:\R\mapsto\R^+\cup\{0\}$ satisfies:  $\hbox{supp } \eta\subset(1,2)$ and is normalized to one, namely  $\int_\R \eta(s)  ds=1$. Observe that we have  $\eqref{vir-1}\leq R^{-2}\|\Delta^2_{\bar x}\psi\|_{L^\infty_{\bar x}}M=o_R(1)$, while the \emph{local} term \eqref{vir-2} can be estimated as in Martel's paper \cite{Mar} (see also \cite{ADF, DF, BFG-21}  for similar results on different dispersive models).  Precisely,
\[
\eqref{vir-2}\leq o_R(1)+o_R(1)\|\nabla u\|^{\alpha}_{L^2},
\]
hence, by using the Young's inequality we have 
\begin{equation}\label{eq:vir3-4}
\eqref{vir-1}+\eqref{vir-2}\leq o_R(1)+o_R(1)\|\nabla u\|^{\alpha}_{L^2}\lesssim o_R(1)+o_R(1)\|\nabla u\|^{2}_{L^2}.
\end{equation}
Hereafter, we use the \emph{small $o$} notation to refer to negative powers of $R$. \medskip

\noindent Our main task  is to show that we can handle in a suitable way also the \emph{non-local} contribution $\eqref{vir-3}+\eqref{vir-4}$. By its definition, we get that the function $\psi_R$ fulfils  
\[
\nabla_{\bar x}\psi_R(x)=
\begin{cases}
2\bar x \quad &\hbox{for } \quad |\bar x|^2\leq R^2\\
0 \quad &\hbox{for } \quad  |\bar x|^2>2R^2
\end{cases},
\]
hence $\hbox{supp }\nabla_{\bar x}\psi_R$  is contained in the cylinder of radius $\sqrt 2 R.$
We split the function $u$ by cutting-off it in the interior and in the exterior of  a cylinder of radius $4R$, namely we write  $u=u_{\leq 4R}+u_{\geq 4R}$ where 
\[u_{\geq 4R}=\bold 1_{\{|\bar x|\leq 4R\}}u \quad \hbox{ and } \quad u_{\geq 4R}=\bold 1_{\{|\bar x|\geq 4R\}}u.
\] 
Since $\hbox{supp } \nabla \psi_R\cap \hbox{supp }u_{\geq 4R}=\emptyset$ we get
\begin{align}\notag
\int\nabla_{\bar x}\psi_R\cdot\nabla_{\bar x}&\left(  E_1(|u|^2)\right)|u|^2 dx\\\label{eq:oi}
&=\int\nabla_{\bar x}\psi_R\cdot\nabla_{\bar x}\left(  E_1(|u_{\geq 4R}|^2)\right)|u_{\leq 4R}|^2dx\\\label{eq:ii}
&+\int\nabla_{\bar x}\psi_R\cdot\nabla_{\bar x}\left(  E_1(|u_{\leq 4R}|^2)\right)|u_{\leq 4R}|^2 dx.
\end{align}
By integration by parts,
\begin{align}\label{eq:oi2}
\eqref{eq:oi}&=-\int\Delta_{\bar x}\psi_R  E_1(|u_{\geq 4R}|^2)|u_{\leq 4R}|^2 dx\\\label{eq:oi3}
&-\int\nabla_{\bar x}\psi_R\cdot\nabla_{\bar x}\left(|u_{\leq 4R}|^2\right)  E_1(|u_{\geq 4R}|^2) dx;
\end{align}
by noting that $\|\Delta_{\bar x}\psi_R\|_{L^\infty}\lesssim1$ and by using \eqref{eq:est-r4}, we obtain
\begin{equation}\label{eq:oi4}
\eqref{eq:oi2}\lesssim R^{-3}\|u_{\geq 4R}\|_{L^2}^2\|u_{\leq 4R}\|_{L^2}^2\lesssim R^{-3}M^2\lesssim R^{-3}.
\end{equation}
Similarly, this time by using that  $|\nabla_{\bar x}\psi_R|\lesssim R$ on its support, \eqref{eq:est-r4} and the Cauchy-Schwarz's inequality give
\begin{equation}\label{eq:oi5}
\eqref{eq:oi3}\lesssim  R^{-2}\|u_{\geq 4R}\|_{L^2}^2\|u_{\leq 4R}\|_{L^2}\|u_{\leq 4R}\|_{\dot H^1}\lesssim R^{-2}\|u(t)\|_{\dot H^1}.
\end{equation}
By \eqref{eq:oi2},  \eqref{eq:oi3},  \eqref{eq:oi4},  \eqref{eq:oi5} we get, with the Young's inequality,
\begin{equation}\label{est:oi7}
 \eqref{eq:oi} \lesssim R^{-2}\|u\|_{\dot H^1}+R^{-3}\lesssim o_R(1)+o_R(1)\|u(t)\|_{\dot H^1}^2.
\end{equation}

\noindent We move to the estimate for \eqref{eq:ii}. By setting $\tilde\psi_R=\psi_R-|\bar x|^2$ we rewrite 
\begin{align}\label{eq:split-V}
\eqref{eq:ii}&
=\int\nabla_{\bar x}\tilde\psi_R\cdot\nabla_{\bar x}\left(  E_1(|u_{\leq 4R}|^2)\right)|u_{\leq 4R}|^2 dx\\\label{eq:split-V2}
&+2\int\bar x\cdot\nabla_{\bar x}\left(  E_1(|u_{\leq 4R}|^2)\right)|u_{\leq 4R}|^2 dx.
\end{align}
We further localize the function $u_{\leq 4R}$ by splitting  $u_{\leq 4R}=u_{\leq R/10}+u_{R/10<4R},$ where 
\[
u_{\leq R/10}=\bold 1_{\{|\bar x|\leq R/10\}}u \quad \hbox{ and } \quad u_{R/10<4R}=\bold 1_{\{R/10\leq|\bar x|\leq 4R\}}u.
\]
Note that $\hbox{supp }\nabla_{\bar x}\tilde\psi_R\subset \{|\bar x|\geq R\},$ hence $\hbox{supp }\nabla_{\bar x}\tilde\psi_R\cap \{|\bar x|\leq R/10\}=\emptyset$. Therefore we can write
\begin{align}\label{eq:split-V3}
\hbox{R.H.S.}\eqref{eq:split-V}&=\int\nabla_{\bar x}\tilde\psi_R\cdot\nabla_{\bar x}\left(  E_1(|u_{\leq R/10}|^2)\right)|u_{R/10<4R}|^2 dx\\\label{eq:split-V4}
&+\int\nabla_{\bar x}\tilde\psi_R\cdot\nabla_{\bar x}\left(  E_1(|u_{R/10<4R}|^2)\right)|u_{R/10<4R}|^2 dx.
\end{align}
After an integration by parts, the R.H.S. of \eqref{eq:split-V3} is controlled as \eqref{eq:oi} (see \eqref{est:oi7}):
\begin{equation*}
\hbox{R.H.S.}\eqref{eq:split-V3}
\lesssim o_R(1)+o_R(1)\|u(t)\|_{\dot H^1}^2.
\end{equation*}

\noindent It remains to prove a suitable estimate for the term \eqref{eq:split-V4}.
By setting $g=|u_{R/10<4R}|^2$ and by making use of the Plancherel identity we get, with $\bar \xi=(\xi_2,\xi_3)$, 
\begin{align}\notag
\eqref{eq:split-V4}&
=\iint \hat g(\eta)\widehat{\nabla_{\bar x}\tilde\rho_R}(\xi-\eta)\cdot\left(\frac{\xi_j\bar\xi}{|\xi|}+\frac{\eta_j\bar\eta}{|\eta|}-\frac{\eta_j\bar\eta}{|\eta|}\right)\frac{\xi_j}{|\xi|}\hat g(\xi)\,d\eta\,d\xi\\\label{eq:3.8}
&=-\frac12\int\Delta_{\bar x}\tilde\psi_R|  E_1 g(x)|^2 dx\\\label{eq:3.9}
&+\iint \hat g(\eta)\widehat{\nabla_{\bar x}\tilde\psi_R}(\xi-\eta)\cdot\left(\frac{\xi_1\bar\xi}{|\xi|}-\frac{\eta_1\bar\eta}{|\eta|}\right)\frac{\xi_1}{|\xi|}\hat g(\xi)  d\eta  d\xi.
\end{align}
As $\|\Delta_{\bar x}\tilde\psi_R\|_{L^\infty}\lesssim 1$, the $L^2\mapsto L^2$ continuity of $E_1$ gives: 
\begin{equation}\label{eq:sec5:B-1}
\eqref{eq:3.8}\lesssim \|u\|^4_{L^4(|\bar x|\geq R/10)}\lesssim R^{-1}\|u\|_{\dot H^1}^2,
\end{equation}
where we used  \eqref{eq:strauss} in the final step.   As for the term  \eqref{eq:3.9}, we explicitly compute 
\begin{align*}
\eqref{eq:3.9}&
=\int \frac{\xi_1}{|\xi|}\hat g(\xi)\int\hat g(\xi_1, \bar \eta)\widehat{\tilde\psi_R}(\bar\xi-\bar\eta)(\bar\xi-\bar\eta)\cdot\left(\frac{\xi_1\bar\xi}{|\xi|}-\frac{\eta_1\bar\eta}{\sqrt{\xi_1^2+|\bar\eta|^2}}\right)  d\bar\eta  d\xi. 
\end{align*}
Note that by the mean value theorem, $\left|\frac{\xi_1\bar\xi}{|\xi|}-\frac{\eta_1\bar\eta}{\sqrt{\xi_1^2+|\bar\eta|^2}}\right|\lesssim |\bar \xi-\bar\eta|$, hence 
\begin{align}\notag
\eqref{eq:3.9}&\lesssim \int|\hat g(\xi)|\int|\hat g(\xi_1, \bar \eta)|\left|\widehat{\tilde\psi_R}(\bar\xi-\bar\eta)\right||\bar\xi-\bar\eta|^2  d\bar\eta  d\xi\\\notag
&\leq\int|\hat g(\xi)|\int|\hat g(\xi_1, \bar \eta)|\left|\widehat{\Delta_{\bar x}\psi_R}(\bar\xi-\bar\eta)\right|  d\bar\eta  d\xi\\\notag
&+\int|\hat g(\xi)|\int|\hat g(\xi_1, \bar \eta)|\left|\widehat{\Delta_{\bar x}|\bar x|^2}(\bar\xi-\bar\eta)\right|  d\bar\eta  d\xi\\\label{eq:fin-four}
&=\int |\hat g(\xi)|\left( |\hat g( \xi_1, \cdot)|\ast \left|\widehat{\Delta_{\bar x}\psi_R}\right|\right)(\bar \xi)  d\xi\\\label{eq:fin-four2}
&+4\int |\hat g(\xi)|^2  d\xi.
\end{align}
By the definition of the functions $g$ and $u_{R/10<4R}$, and by the isometry property of the Fourier transform, we easily bound, again by using \eqref{eq:strauss},
\begin{equation}\label{eq:delta}
\eqref{eq:fin-four2}\lesssim \|u\|^4_{L^4(|\bar x|\geq R/10)}\lesssim R^{-1}\|u\|_{\dot H^1}^2.
\end{equation}
As for the remaining term \eqref{eq:fin-four}, we have: 
\begin{equation}\label{eq:almdo}
\begin{aligned}
\eqref{eq:fin-four}&\leq\int_{\R}\|\hat g(\xi_1,\cdot)\|_{L^2
(\R^2_{\bar\xi})}\left\| |\hat g( \xi_1, \cdot)|\ast \left|\widehat{\Delta_{\bar x}\psi_R}\right| \right\|_{L^2(\R^2_{\bar\xi})}  d\xi_1\\
&\leq \|\widehat{\Delta_{\bar x}\psi_R}\|_{L^1(\R^2_{\bar\xi})}\int_{\R}\|\hat g( \xi_1, \cdot)\|_{L^2(\R^2_{\bar\xi})}^2  d\xi_1= \|\widehat{\Delta_{\bar x}\psi_R}\|_{L^1(\R^2_{\bar\xi})}\|\hat g\|_{L^2}^2\\
&=\|\widehat{\Delta_{\bar x}\psi_R}\|_{L^1(\R^2_{\bar\xi})}\| g\|_{L^2}^2\lesssim \|u\|_{L^4(|\bar x|\geq R/10)}^4\lesssim R^{-1}\|u\|_{\dot H^1}^2.
\end{aligned}
\end{equation}
where in order we used: the Cauchy-Schwarz's inequality and the Young's inequality for convolutions with respect to $\bar \xi$, the Cauchy-Schwarz's inequality with respect to $\xi_1$,  the Fourier $L^2$ isometry property, and the fact that  $\widehat{\Delta_{\bar x}\psi_R}$ is integrable (with bound independent of $R$), as it is the Fourier  transform of a compactly supported regular function. Again, the norm of $u$ outside a cylinder is estimated by \eqref{eq:strauss}. 
By glueing the estimates above, we see that the \emph{non-local} term $\eqref{vir-3}+\eqref{vir-4}$ is estimated by
\begin{align*}
&\eqref{vir-3}+\eqref{vir-4}\leq o_R(1)+o_R(1)\|u\|_{\dot H^1}^2\\
&\quad+4c_2\left(\int\bar x\cdot\nabla_{\bar x}\left(  E_1(|u_{\leq 4R}|^2)\right)|u_{\leq 4R}|^2 dx+\int x_1\partial_{x_1}\left(  E_1(|u|^2)\right)|u|^2 dx\right).
\end{align*}
A straightforward computation gives 
\begin{align*}
&\int\bar x\cdot\nabla_{\bar x}\left(  E_1(|u_{\leq 4R}|^2)\right)|u_{\leq 4R}|^2 dx+\int x_1\partial_{x_1}\left(  E_1(|u|^2)\right)|u|^2 dx\\
&=\int x\cdot\nabla\left(  E_1(|u_{\leq 4R}|^2)\right)|u_{\leq 4R}|^2 dx+\int x_1\partial_{x_1}\left(  E_1(|u_{\leq 4R}|^2)\right)|u_{\geq 4R}|^2 dx\\
&+\int x_1\partial_{x_1}\left(  E_1(|u_{\geq 4R}|^2)\right)|u_{\leq 4R}|^2 dx+\int x_1\partial_{x_1}\left(  E_1(|u_{\geq 4R}|^2)\right)|u_{\geq 4R}|^2 dx\\
&=-\frac32\int \left(  E_1(|u_{\leq 4R}|^2)\right)|u_{\leq 4R}|^2 dx+\int x_1\partial_{x_1}\left(  E_1(|u_{\leq 4R}|^2)\right)|u_{\geq 4R}|^2 dx\\
&+\int x_1\partial_{x_1}\left(  E_1(|u_{\geq 4R}|^2)\right)|u_{\leq 4R}|^2 dx+\int x_1\partial_{x_1}\left(  E_1(|u_{\geq 4R}|^2)\right)|u_{\geq 4R}|^2 dx,
\end{align*}
where we used the identity $2\int x\cdot\nabla\left(  E_1(f)\right)f dx=-3\int   E_1(f) f dx$, see \cite{Cipo}. We estimate now the  four terms above. By using twice  the Plancherel theorem, we can compute 
\begin{align}\notag
\int x_1\partial_{x_1}&\left(  E_1(|u_{\geq 4R}|^2)\right)|u_{\geq 4R}|^2 dx\\\label{eq:est-oo}
&=-\frac12\int  E_1(|u_{\geq 4R}|^2)|u_{\geq 4R}|^2 dx\\\label{eq:est-oo-2}
&-\frac12\int\xi_1\partial_{\xi_1}\left(\frac{\xi_1^2}{|\xi|^2}\right)\widehat{ |u_{\geq 4R}|^2}\overline{ \widehat {|u_{\geq 4R}|^2}}  d\xi,
\end{align}
while, by similar computations, once passed in the frequency space, we get 
\begin{align}\notag
\int x_1\partial_{x_1}\left(  E_1(|u_{\leq 4R}|^2)\right)&|u_{\geq 4R}|^2 dx+\int x_1\partial_{x_1}\left(  E_1(|u_{\geq 4R}|^2)\right)|u_{\leq 4R}|^2 dx\\\label{eq:term:oo}
&=-\int\left(  E_1(|u_{\leq 4R}|^2)\right)|u_{\geq 4R}|^2 dx\\\label{eq:term:oo-bis}
&-\int\xi_1\left(\partial_{\xi_1}\frac{\xi_1^2}{|\xi|^2}\right)\widehat{ |u_{\leq 4R}|^2 }\overline{\widehat{|u_{\geq 4R}|^2}}  d\xi.
\end{align}
We explicitly write $\xi_1\partial_{\xi_1}\left(\frac{\xi_1^2}{|\xi|^2}\right)$ and we observe that it is bounded:
\begin{equation}\label{eq:der-ker-bound1}
\xi_1\partial_{\xi_1}\left(\frac{\xi_1^2}{|\xi|^2}\right)=\frac{2\xi_1^2(\xi_2^2+\xi_3^2)}{|\xi|^4}=\frac{2\xi_1^2}{|\xi|^2}-\frac{2\xi_1^4}{|\xi|^4}\leq 4.
\end{equation}
By Remark \ref{rem:bound}, the boundedness of the above Fourier symbol implies an $L^2\mapsto L^2$ continuity, hence  $\eqref{eq:est-oo}+\eqref{eq:est-oo-2}$ is simply estimated, jointly  with \eqref{eq:strauss}, by
\[
\eqref{eq:est-oo}+\eqref{eq:est-oo-2} \lesssim \|u_{\geq 4R}\|_{L^4}^4\lesssim R^{-1}\|u\|_{\dot H^1}^2.
\]
We are left with $\eqref{eq:term:oo}+\eqref{eq:term:oo-bis}$. First of all we note by \eqref{eq:der-ker-bound1} that, 
up to constants, $\xi_1\partial_{\xi_1}\left(\frac{\xi_1^2}{|\xi|^2}\right)$ is the sum of symbols defining the pseudo-differential operators $E_1$ and $E_1^2$.  This in turn implies that $\eqref{eq:term:oo}+\eqref{eq:term:oo-bis}$ can be rewritten as 
\[
-3\int\left(  E_1( |u_{\leq 4R}|^2)\right)|u_{\geq 4R}|^2 dx+2\int \left(  E_1^2(|u_{\leq 4R}|^2)\right)|u_{\geq 4R}|^2 dx.
\]
Again by splitting $u_{\leq 4R}=u_{\leq R/10}+u_{R/10<4R}$ 
we decompose  
\begin{equation*}
\begin{aligned}
\int\left(  E_1( |u_{\leq 4R}|^2)\right)|u_{\geq 4R}|^2 dx&=\int\left(  E_1( |u_{R/10<4R}|^2)\right)|u_{\geq 4R}|^2 dx\\
&+\int\left(  E_1(|u_{\leq R/10}|^2)\right)|u_{\geq 4R}|^2 dx,
\end{aligned}
\end{equation*}
and by using the Cauchy-Schwarz's inequality, the $L^2\mapsto L^2$ continuity of  $E_1$, and \eqref{eq:strauss},  we obtain
\begin{equation}\label{eq:wio-uo}
\int  E_1( |u_{R/10<4R}|^2)|u_{\geq 4R}|^2 dx\lesssim \|u\|_{L^4(|\bar x|\geq R/10)}^4\lesssim R^{-1}\|u\|_{\dot H^1}^2.
\end{equation}
By using  \eqref{eq:est-r4} we instead give the bound
\begin{equation*}
\int  E_1( |u_{\leq R/10}|^2)|u_{\geq 4R}|^2 dx\lesssim R^{-3}\|u_{\leq R/10}\|_{L^2}^2\|u_{\geq 4R}\|^2_{L^2}\lesssim R^{-3}.
\end{equation*}
With the same decomposition of the function $u_{\leq 4R}$, we separate the term defined by $E_1^2$ as 
\begin{equation*}
\begin{aligned}
\int   E_1^2(|u_{\leq 4R}|^2)|u_{\geq 4R}|^2 dx&=\int   E_1^2(|u_{\leq R/10}|^2)|u_{\geq 4R}|^2 dx\\
&+\int   E_1^2(|u_{R/10<4R}|^2)|u_{\geq 4R}|^2 dx;
\end{aligned}
\end{equation*}
then, similarly to \eqref{eq:wio-uo}, we  have
\begin{equation*}
\int \left(  E_1^2(|u_{R/10<4R}|^2)\right)|u_{\geq 4R}|^2 dx\lesssim R^{-1}\|u\|_{\dot H^1}^2,
\end{equation*}
while, by using \eqref{eq:est-r4}, we can control  
\begin{equation*}
\int   E_1^2(|u_{\leq R/10}|^2)|u_{\geq 4R}|^2 dx\lesssim R^{-3}\|u_{\leq R/10}\|_{L^2}^2\|u_{\geq 4R}\|_{L^2}^2\lesssim R^{-3}.
\end{equation*}
The final term to deal with is $-3\int \left(  E_1(|u_{\leq 4R}|^2)\right)|u_{\leq 4R}|^2 dx$. By adding and subtracting $|u|^2$ to $|u_{\leq 4R}|^2$, we can fall back into  the same localized objects as in the above discussions, hence we get
\begin{equation}\label{eq:very-last}
\begin{aligned}
-\frac32\int \left(  E_1(|u_{\leq 4R}|^2)\right)|u_{\leq 4R}|^2 dx&=-\frac32\int   E_1(|u|^2)|u|^2 dx\\
&+ o_R(1)+o_R(1)\|u\|_{\dot H^1}^2.
\end{aligned}
\end{equation}
At this point, by collecting the above estimates, we end-up with
\begin{equation}\label{eq:35-36}
\eqref{vir-3}+\eqref{vir-4}\leq o_R(1)+o_R(1)\|u\|^2_{\dot H^1}-6c_2\int E_1(|u|^2)|u|^2 dx.
\end{equation}
We can now summarize all the previous contributions towards the conclusion of the proof. From  \eqref{vir0},\eqref{vir-1},\eqref{vir-2},\eqref{vir-3}, and \eqref{vir-4}, coupled with \eqref{eq:vir3-4} and \eqref{eq:35-36}, and by recalling the definition of $P$, see \eqref{eq:G}, we have 
\begin{equation*}
\begin{aligned}
\frac{d^2}{dt^2}V_{x_1^2+\psi_R(\bar x)}(t)&\leq 8\left(\int|\nabla u(t)|^2 dx -\frac{3c_1\alpha}{2(\alpha+2)}\int |u(t)|^{\alpha+2} dx\right) \\
&-6c_2\int E_1(|u(t)|^2)|u(t)|^2 dx +o_R(1)+o_R(1)\|u(t)\|^2_{\dot H^1}\\
&=8P(u(t)) +o_R(1)+o_R(1)\|u(t)\|^2_{\dot H^1}.
\end{aligned}
\end{equation*}
Note that from Proposition \ref{prop:P-nega} and the Sobolev embedding, it can be claimed that $\inf_{(-T_{-}, T_{+})}\|u(t)\|_{\dot H^1}\geq\beta>0$: otherwise, by contradiction, along a sequence of times $\{t_n\}\subset(-T_{-}, T_{+})$ we would have that $P(u(t_n))\to0$. Thus, provided we chose $R\gg1$, from the estimates above and Proposition \ref{prop:P-grad}, we get that $\frac{d^2}{dt^2}V_{x_1^2+\psi_R(\bar x)}(t)\lesssim-1$. A  convexity argument concludes the proof.


\section*{Acknowledgements}
\noindent  The author was supported by the EPSRC New Investigator Award (grant no. EP/S033157/1).


\end{document}